\documentclass[11pt,a4paper]{article}
\usepackage[utf8]{inputenc}
\usepackage[english]{babel}
\usepackage{amsmath,amssymb,amsfonts,amsthm}
\usepackage{geometry}
\usepackage{graphicx}
\usepackage{microtype}
\usepackage{booktabs}
\usepackage{cite}
\usepackage{hyperref}

\newtheorem{theorem}{Theorem}[section]
\newtheorem{corollary}[theorem]{Corollary}
\newtheorem{lemma}[theorem]{Lemma}
\theoremstyle{definition}
\newtheorem{remark}[theorem]{Remark}
\newtheorem{example}[theorem]{Example}

\numberwithin{equation}{section}

\begin{document}

\title{\textbf{\Large A sufficient condition for the Liouville property of solutions to the Euler-Lagrange equations of variational integrals with nonstandard growth in terms of the ellipticity modulus near the origin}}
\author{\textbf{Michael Bildhauer} \quad \textbf{Martin Fuchs}}
\date{}
\maketitle

\begin{abstract}
\noindent We consider entire solutions $u$: $\mathbb{R}^{n}\rightarrow\mathbb{R}$ of the Euler-Lagrange equation $\operatorname{div}(\nabla f(\nabla u))=0$ 
with convex density $f$: $\mathbb{R}^{n}\rightarrow\mathbb{R}$ of nonstandard growth and nonuniform ellipticity. 
It turns out that the Liouville property can be characterized through the behavior of $D^{2}f(p)$ for points $p\in\mathbb{R}^{n}$ close to the origin:
if, for example, $D^{2}f(p)$ behaves as $|p|^{\rho-2}\mathbf{1}$, $\mathbf{1}=\text{unit matrix}$, with exponent $\rho\ge2$ and for $|p|\le1$ then the validity of 
$n\le\rho+2$ together with the boundedness of $u$ implies its constancy.

\vspace{3mm}
\end{abstract}

\section{Introduction}

In this note we discuss the Liouville property for entire solutions $u\in C^{2}(\mathbb{R}^{n})$ of the nonlinear elliptic equation
\begin{equation}
\operatorname{div} (\nabla f(\nabla u))=0
\end{equation}
arising as the Euler-Lagrange equation of the variational integral $\int f(\nabla u) \, dx$ with convex density 
$f:\mathbb{R}^{n}\rightarrow\mathbb{R}$ of class $C^{2}$. Here we mean by Liouville property the fact that any bounded entire solution 
$u$: $\mathbb{R}^{n}\rightarrow\mathbb{R}$ of (1.1) must be a constant function.\\

Surprisingly, it turns out that a positive answer just depends on the behavior of $D^{2}f(p)$ for $p$ close to the origin, for example, if $D^{2}f(p)$ is of order $|p|^{\rho-2}$ as 
$p\rightarrow 0$ with exponent $\rho\ge2$, then we have the Liouville property provided that $n\le\rho+2$, and this is true under rather general hypotheses 
concerning the growth of $f$ at infinity.\\ 

Before going into details we like to remark that the study of equations with nonstandard growth and nonuniform ellipticity is a very active field, 
which has attracted the research of many scientists over the last decades focussing on the questions of existence and regularity of solutions. 
For a detailed overview including historical comments the reader is referred to [1]. However, we could not find a discussion of the Liouville property 
in this general setting, but it should be mentioned that there are contributions concerning quasilinear or uniformly elliptic equations, 
which can be found for instance in [2], [3], [4], [5] and in the references quoted therein. We also advise the reader to consult the textbooks [6], [7].\\

Let us now state our assumptions concerning the density $f\in C^{2}(\mathbb{R}^{n})$ occurring in equation (1.1): we have
\begin{equation}
D^{2}f(p)(q,q)>0, \quad p,q\in\mathbb{R}^{n}, \, q\ne0, \quad \nabla f(0)=0,
\end{equation}
implying the strict convexity of $f$. We further require ($c,c_{1},\dots$ denoting positive constants)
\begin{equation}
|\nabla f(p)|\le c|p|^{m-1} \quad \text{for } p\in\mathbb{R}^{n}, \, |p|\ge1
\end{equation}
with an arbitrary exponent $m>1$, and addition we assume
\begin{equation}
p\cdot\nabla f(p)\ge c_{1}|p|^{-\mu}|p||\nabla f(p)|, \quad p\in\mathbb{R}^{n}, \, |p|\ge1
\end{equation}
for some number $\mu\in[0,1)$ for which we also have
\begin{equation}
|D^{2}f(p)|\ln^{2}(1+|p|)\le c_{2}|p|^{-\mu}|\nabla f(p)|, \quad p\in\mathbb{R}^{n}, \, |p|\ge1.
\end{equation}

The conditions (1.3)--(1.5) are conditions at infinity, since the bound $|p|\ge1$ for $p\in\mathbb{R}^{n}$ can be replaced by $|p|\ge M$ for some arbitrary number $M$.
 Inequality (1.3) excludes exponential growth, but the exponent $m$ may be as large as we want. 
 Condition (1.4) states that the radial part of $\nabla f(p)$, i.e. the quantity $\nabla f(p)\cdot\frac{p}{|p|}$, is not allowed to degenerate at infinity 
 in the sense that we have $|p|^{-\mu}|\nabla f(p)|$ as a lower bound for some exponent $\mu\in[0,1)$. 
 Finally, the growth of $|D^{2}f(p)|$ as $p\rightarrow\infty$ is limited in terms of $|\nabla f(p)|$ equipped 
 with the weight function $|p|^{-\mu}$, which is made precise in (1.5). We refer to Section 2 for a list of prominent examples.\\

The conditions (1.3)--(1.5) are completed by
\begin{equation}
\exists\rho\ge2: \quad c_{3}|p|^{\rho-2}|q|^{2}\le D^{2}f(p)(q,q)\le c_{4}|p|^{\rho-2}|q|^{2}, \quad p,q\in\mathbb{R}^{n}, \, |p|\le1
\end{equation}
and we can state\\

\begin{theorem} \label{thm1.1}
Let (1.2)--(1.6) hold and let $u\in C^{2}(\mathbb{R}^{n})$ denote an entire solution of equation (1.1) for which
\begin{equation}
\lim_{x\rightarrow\infty}\frac{|u(x)|}{|x|^{\beta}}<\infty
\end{equation}
with
\begin{equation}
\beta:=\beta(\rho):=\frac{1}{\rho}(\rho+2-n).
\end{equation}
Then $u$ is an affine function, i.e. $\partial_{i}\partial_{j}u=0$ on $\mathbb{R}^{n}$ for $i, j=1,\dots,n.$
\end{theorem}

If $n=2$, then $\beta=1$ and in conclusion:\\

\begin{corollary}
Under the assumptions of Theorem \ref{thm1.1} consider the case $n=2$ and let $u \in C^{2}(\mathbb{R}^{2})$ denote an entire solution of (1.1) such that
\begin{equation*}
\lim_{x\rightarrow\infty}\frac{|u(x)|}{|x|}<\infty.
\end{equation*}  
Then $u$ is an affine function.\\
\end{corollary}

Note that the number $\beta=\beta(\rho)$ from (1.8) is $\ge 0$ if and only if $n\le\rho+2$. This implies the Liouville property for dimensions 
$n\le n(\rho)$ with upper bound $n(\rho):=2+\rho$ determined just through the behavior of $D^{2}f(p)$ for $p$ close to the origin as expressed in (1.6).\\

\begin{corollary}[Liouville property]
Let (1.2)--(1.6) hold and assume that $n\le\rho+2$. Then the constant functions are the only bounded entire solutions of equation (1.1).\\
\end{corollary}

In fact, if $u\in L^{\infty}(\mathbb{R}^{n})$ and $n\le\rho+2$, then (1.7) clearly holds, hence $u$ is an affine function and thereby constant.\\

\begin{remark}
It is possible to replace (1.6) by the anisotropic condition
\begin{equation*}
c_{5}|p|^{t-2}|q|^{2}\le D^{2}f(p)(q,q)\le c_{6}|p|^{\rho-2}|q|^{2}, \quad p,q\in\mathbb{R}^{n}, \, |p|\le1,
\end{equation*}
with exponents $2\le\rho\le t$. Then Theorem \ref{thm1.1} continues to hold, provided we define in (1.8)
\begin{equation*}
\beta:=\beta(t,\rho):=\frac{1}{T}(2+T-n), \quad T:=\frac{t}{t-\rho+1}
\end{equation*}
as the appropriate value of $\beta$, and if we require $\rho+2\ge t.$ Some comments will be given in Section 4.
\end{remark}

\begin{remark}
One might wonder about the occurrence of the quantity $\ln^{2}(1+|p|)$ on the left-hand side of inequality (1.5). 
Of course we can replace (1.5) by the stronger condition (keeping (1.4)) that
\begin{equation*}
|D^{2}f(p)|\le c_{2}|p|^{-(\mu+\epsilon)}|\nabla f(p)|, \quad p\in\mathbb{R}^{n}, \, |p|\ge1
\end{equation*}
holds with exponent $\mu$ from (1.4) and for some arbitrarily small number $\epsilon>0$. However, we prefer to work with a less restrictive assumption, 
and at the same time it turns out that the expression $|D^{2}f(\nabla u)|\ln^{2}(1+|\nabla u|)$ 
occurs in a limit version of a Caccioppoli type inequality, which will be presented in Lemma 3.1 being valid for solutions $u$ 
of equation (1.1) under very general hypotheses concerning the density $f$.\\
\end{remark}

Our paper is organized as follows: in Section 2 we collect some examples to which Theorem \ref{thm1.1} and the corollaries apply. 
Section 3 is devoted to an inequality of Caccioppoli type, which is crucial for the proof of Theorem \ref{thm1.1} presented in Section 4. 
The Caccioppoli inequality Lemma 3.1 will explain the occurrence of the factor $\ln^{2}(1+|p|)$ on the left-hand side of (1.5) and 
during the proof of Lemma 3.1 the comments from Remark 1.5 will become more transparent.

\section{Examples}

We present some examples of densities $f$: $\mathbb{R}^{n}\rightarrow\mathbb{R}$ to which Theorem \ref{thm1.1} can be applied.\\

\begin{example}[anisotropic power growth]
Let $f\in C^{2}(\mathbb{R}^{n})$ satisfy (1.2) and (1.6). In addition we require: there exist constants $c_{1}$, $c_{2}$ and exponents $1<\tau\le r$ such that
\begin{equation}
c_{1}|p|^{\tau-2}|q|^{2}\le D^{2}f(p)(q,q)\le c_{2}|p|^{r-2}|q|^{2}
\end{equation}
holds for $p, q\in\mathbb{R}^{n}$, $|p|\ge1$. Then we claim: if
\begin{equation}
r<\tau+\frac{1}{2},
\end{equation}
then we have (1.3), (1.4), (1.5) for some $\mu\in[0,1)$ and with arbitrary $m\ge r$.\\ 

The validity of (1.3) for $m\ge r$ is immediate by the second inequality in (2.1). For (1.4) and (1.5) we choose $p\in\mathbb{R}^{n}$, $|p|\ge1$, and observe that by (2.1)
\begin{align*}
|\nabla f(p)||p| &\ge \nabla f(p)\cdot p=\int_{0}^{1}D^{2}f(tp)(p,p)dt \\
&\ge c_{1}\int_{1/2}^{1}t^{\tau-2}dt|p|^{\tau}
\end{align*}
for $|p|\ge2$, whereas in the case $1\le|p|\le2$ we use (1.6) to get
\begin{equation}
\begin{cases}
\nabla f(p) \cdot p \ge c|p|^\tau, \quad |\nabla f(p)| \ge c |p|^{\tau-1}, \\[2ex]
p \in \mathbb{R}^n, \; |p| \ge 1,
\end{cases}
\end{equation}
with a positive constant $c$. The second part of (2.1) implies with the help of (2.3)
\begin{align*}
|D^{2}f(p)|\ln^{2}(1+|p|) &\le c \ln^{2}(1+|p|)|p|^{r-2} \\[2ex]
&=c \ln^{2}(1+|p|)|p|^{\tau-1}|p|^{r-\tau-1} \\[2ex]
&\le c \ln^{2}(1+|p|)|p|^{r-\tau-1}|\nabla f(p)| \\[2ex]
&\le c|p|^{-\mu_{1}}|\nabla f(p)|
\end{align*}  
for any $\mu_{1} \in [0, \tau+1-r)$.\\

Moreover, again for $|p|\ge1$, quoting (1.6) and (2.1)
\begin{align*}
|\nabla f(p)| &= \left|\int_{0}^{1}\frac{d}{dt}\nabla f(tp)dt\right| \le \int_{0}^{1}|D^{2}f(tp)|\,dt|p| \\[2ex]
&= |p|\left(\int_{0}^{1/|p|}|D^{2}f(tp)|dt+\int_{1/|p|}^{1}|D^{2}f(tp)|dt\right) \\[2ex]
&\le c|p|\left(\frac{1}{|p|}+|p|^{r-2}\right)
\end{align*}
so that $|\nabla f(p)|\le c|p|^{r-1}$ and in conclusion (see (2.3))
\begin{align*}
p\cdot\nabla f(p) &\ge c|p|^{\tau}=c|p|^{r-1}|p||p|^{\tau-r+1} \\[2ex]
&\ge c|p|^{\tau-r}|\nabla f(p)||p| \ge c|p|^{-\mu_{2}}|\nabla f(p)||p|
\end{align*}
provided that $\mu_{2}\ge r-\tau.$\\ 

In order to find a common value $\mu$ for $\mu_{1}$ and $\mu_{2}$ we need to know that $r-\tau<\tau+1-r$, 
which is just a consequence of (2.2) completing the discussion of (1.4) and (1.5).\\
\end{example}

\begin{example}[densities $f$ with linear growth]
Again we consider $f\in C^{2}(\mathbb{R}^{n})$ with (1.2) and (1.6). We further assume that
\begin{equation}
a(|p|-1)\le f(p)\le A(|p|+1) \quad p\in\mathbb{R}^{n}
\end{equation}
holds with constants $a, A>0$. Moreover, we require
\begin{equation}
D^{2}f(p)(q,q)\le\Theta(|p|)|q|^{2} \quad p,q\in\mathbb{R}^{n}, \, |p|\ge1,
\end{equation}
for a function $\Theta:[0,\infty)\rightarrow[0,\infty)$ such that
\begin{equation}
\limsup_{t\rightarrow\infty}\Theta(t)\ln^{2}(t)<\infty.
\end{equation}

Then we are in the situation of Theorem \ref{thm1.1}, i.e. the conditions (1.3)--(1.6) hold: 
for (1.3) we just observe that from the convexity of $f$ (see (1.2)) together with (2.4) we immediately get the boundedness of $\nabla f$.\\ 

Let
\begin{equation*}
\vartheta(r)=\inf_{|p|=r}\left(\inf_{|q|=1}D^{2}f(p)(q,q)\right), \quad \vartheta:[0,\infty)\rightarrow[0,\infty).
\end{equation*} 
Again from (1.2) we deduce that $\vartheta(r)>0$ for $r>0$. Observing this, it follows for $|p|\ge1$:
\begin{align*}
|p||\nabla f(p)| &\ge p\cdot\nabla f(p)=\int_{0}^{1}D^{2}f(tp)(p,p)dt \\[2ex]
&\ge|p|^{2}\int_{0}^{1}\vartheta(t|p|)dt=|p|\int_{0}^{|p|}\vartheta(s)ds \\[2ex]
&\ge|p|\int_{0}^{1}\vartheta(s)ds =: \alpha|p|
\end{align*}
for some positive constant $\alpha$, and the boundedness of $\nabla f$ clearly implies (1.4) with the choice $\mu=0$. Condition (1.5) (with $\mu=0$) is a consequence of (2.5) and (2.6).\\
\end{example}

\begin{example}[densities with radial structure]
We let $f(p):=g(|p|)$, $p\in\mathbb{R}^{n}$ for a function $g\in C^{2}([0,\infty))$ such that
\begin{equation}
\left\{
\begin{aligned}
&g^{\prime}(0)=0, \\[2ex]
&g^{\prime\prime}(t)>0 \quad \text{for } t>0
\end{aligned}
\right\},
\quad c_{1}t^{\rho-2}\le g^{\prime\prime}(t)\le c_{2}t^{\rho-2} ,\quad 0\le t\le1,
\end{equation}
\begin{equation}
\begin{gathered}
\begin{cases}
g^{\prime}(t)\le c_{3}t^{m-1},   t\ge1, \\[2ex]
g^{\prime\prime}(t)\ln^{2}(1+t)\le c_{4}g^{\prime}(t) \text{ for } t\ge1
\end{cases} \\
,
\end{gathered}
\end{equation}
with constants $c_{1},\dots,c_{4}>0$ and exponents $m>1$, $\rho\ge2$.\\ 
We claim that $f$ satisfies the hypotheses of Theorem \ref{thm1.1}: noting that
\begin{eqnarray*}
\nabla f(p) &=& g^{\prime}(|p|)\frac{p}{|p|}, \quad |\nabla f(p)|=g^{\prime}(|p|) \, , \\[2ex] 
D^2 f(p)(q,q) &=& \frac{g'(|p|)}{|p|}\left[|q|^2 - \left(\frac{p}{|p|}\cdot q\right)^2\right] + g''(|p|)\left(\frac{p}{|p|}\cdot q\right)^2
\end{eqnarray*}
the validity of (1.2) and (1.6) is an immediate consequence of (2.7), whereas (1.3) follows from (2.8).\\ 
Let $p\in\mathbb{R}^{n}$, $|p|\ge1.$ Then we get (again quoting (2.8)):
\begin{align*}
|D^{2}f(p)|\ln^{2}(1+|p|) &\le \max\left\{\frac{g^{\prime}(|p|)}{|p|},g^{\prime\prime}(|p|)\right\}\ln^{2}(1+|p|) \\[2ex]
&\le c \max\left\{\frac{g^{\prime}(|p|)}{|p|}\ln^{2}(1+|p|),g^{\prime}(|p|)\right\} \\[2ex]
&\le cg^{\prime}(|p|)=c|\nabla f(p)|,
\end{align*}
hence (1.5) holds for the choice $\mu=0$. From $p\cdot\nabla f(p)=|p||\nabla f(p)|$ we obtain (1.4).
\end{example}

\section{A Caccioppoli type inequality}

A preliminary version of the following lemma can be found in Proposition 6.1 of [8]. Here we extend Proposition 6.1 in some sense to the limit case $\alpha=-1/2$.

\begin{lemma} \label{lem3.1}
Let $f\in C^{2}(\mathbb{R}^{n})$ denote a convex function, which means that $D^{2}f(p)$ is just positive semidefinite for any $p\in\mathbb{R}^{n}$.\\ 

We define $\lambda$: $[1,\infty)\rightarrow[0,\infty)$, $\Lambda$:$[1,\infty)\rightarrow[0,\infty)$ by
\begin{equation}
\lambda(t):=\begin{cases}t^{-1/2},&t\ge2,\\[2ex] \frac{1}{\sqrt{2}}(t-1),&1\le t\le2\end{cases}, \quad
\Lambda(t):=\begin{cases}t^{-1/2}(t-1)\ln^{2}(t),&t\ge2\\[2ex] \frac{1}{\sqrt{2}}\ln^{2}(t),&1\le t\le2\end{cases}
\end{equation}
Let $\Omega\subset\mathbb{R}^{n}$ denote a domain and consider a solution $u\in C^{2}(\Omega)$ of equation (1.1) now on $\Omega$. \\

Then it holds for any $\eta\in C_{0}^{1}(\Omega)$ and for all indices $i\in\{1,\dots,n\}$ ($\Gamma_{i}:=1+|\partial_{i}u|^{2}$):
\begin{eqnarray}
\lefteqn{\int_{\Omega}D^{2}f(\nabla u)(\nabla\partial_{i}u,\nabla\partial_{i}u)\eta^{2}\lambda(\Gamma_{i})\,dx}\nonumber \\[2ex]
& \le & c\left(\int_{\operatorname{spt}(\nabla\eta)}D^{2}f(\nabla u)(\nabla\partial_{i}u,\nabla\partial_{i}u)\eta^{2}\lambda(\Gamma_{i})\,dx\right)^{1/2} \nonumber\\[2ex]
 && \times \left(\int_{\operatorname{spt}(\nabla\eta)}D^{2}f(\nabla u)(\nabla\eta,\nabla\eta)\Lambda(\Gamma_{i})\,dx\right)^{1/2},
\end{eqnarray}
$c$ denoting a constant independent of $\eta$, $\Omega$, $u$, and in conclusion
\begin{equation}
\int_{\Omega}D^{2}f(\nabla u)(\nabla\partial_{i}u,\nabla\partial_{i}u)\eta^{2}\lambda(\Gamma_{i})\,dx \le c\int_{\Omega}D^{2}f(\nabla u)(\nabla\eta,\nabla\eta)\Lambda(\Gamma_{i})\,dx.
\end{equation}
Moreover, on the right-hand side of (3.3) we may replace $\Lambda(\Gamma_{i})$ by $\Lambda(\Gamma)$, where $\Gamma:=1+|\nabla u|^{2}$.
\end{lemma}

\begin{proof}
The last claim follows from the monotonicity of $\Lambda$. Now we fix $\eta\in C_{0}^{1}(\Omega)$ and $i\in\{1,\dots,n\}$. Let 
\[
\Phi:[1,\infty)\rightarrow[0,\infty)\, , \quad \Phi(t):=t^{-1/2}\ln t\, , 
\]
and observe that $\Psi:=\eta^{2}\partial_{i}u\Phi(\Gamma_{i})$ is admissible in the equation
\begin{equation}
\int_{\Omega}D^{2}f(\nabla u)(\nabla\partial_{i}u,\nabla\Psi)\,dx=0.
\end{equation}
Here equation (3.4) is a direct consequence of the weak form of (1.1). From (3.4) we obtain
\begin{eqnarray*}\lefteqn{
\int_{\Omega}D^{2}f(\nabla u)(\nabla\partial_{i}u,\nabla\partial_{i}u)\eta^{2}[\Phi(\Gamma_{i})+2|\partial_{i}u|^{2}\Phi^{\prime}(\Gamma_{i})]\,dx}\\[2ex] 
&=& -2\int_{\operatorname{spt}(\nabla\eta)}D^{2}f(\nabla u)(\nabla\partial_{i}u,\nabla\eta)\eta\partial_{i}u\Phi(\Gamma_{i})\,dx
\end{eqnarray*}
or equivalently
\begin{eqnarray}\lefteqn{
\int_{\Omega}D^{2}f(\nabla u)(\nabla\partial_{i}u,\nabla\partial_{i}u)\eta^{2}[\Phi(\Gamma_{i})+2(\Gamma_{i}-1)\Phi^{\prime}(\Gamma_{i})]\,dx}\nonumber\\[2ex] 
&=& -2\int_{\operatorname{spt}(\nabla\eta)}D^{2}f(\nabla u)(\nabla\partial_{i}u,\nabla\eta)\partial_{i}u\eta\Phi(\Gamma_{i})\,dx.
\end{eqnarray}
Letting $t:=\Gamma_{i}$ ($\ge1$) for the moment, we get for the expression $[\dots]$ in (3.5), using the definition of $\Phi$,
\begin{equation*}
[\dots]=2(t^{-1/2}-t^{-3/2})+t^{-3/2}\ln t.
\end{equation*}
This implies
\begin{equation}
\Phi(t)+2(t-1)\Phi^{\prime}(t)\ge\begin{cases}t^{-1/2},&t\ge2\\[2ex] \frac{1}{\sqrt{2}}(t-1),&1\le t\le2\end{cases}
\end{equation}
as a lower bound for the term on the left-hand side of (3.5). On the right-hand side we decompose
\begin{equation*}
-2\int_{\operatorname{spt}(\nabla\eta)}\dots = -2\int_{\operatorname{spt}(\nabla\eta)\cap[\Gamma_{i}<2]}\dots -2\int_{\operatorname{spt}(\nabla\eta)\cap[\Gamma_{i}\ge2]}\dots =: (I)+(II)
\end{equation*}
discussing the terms (I) and (II) separately: it holds
\begin{eqnarray*}
|\text{integrand of } (II)| &=& 2|\partial_{i}u||\eta|\Phi(\Gamma_{i})|D^{2}f(\nabla u)(\nabla\partial_{i}u,\nabla\eta)| \\[2ex]
&=& 2(\Gamma_{i}-1)^{1/2}\Gamma_{i}^{-1/2}\ln(\Gamma_{i})|\eta||D^{2}f(\nabla u)(\nabla\partial_{i}u,\nabla\eta)| \\[2ex]
&=& 2\left|D^{2}f(\nabla u)\left(|\eta|\Gamma_{i}^{-1/4}\nabla\partial_{i}u,\,\Gamma_{i}^{-1/4}(\Gamma_{i}-1)^{1/2}\ln(\Gamma_{i})\nabla\eta\right)\right| \\[2ex]
&\le & 2\left(D^{2}f(\nabla u)(\nabla\partial_{i}u,\nabla\partial_{i}u)\Gamma_{i}^{-1/2}\eta^{2}\right)^{1/2}\\[2ex]
&& \times \left(D^{2}f(\nabla u)(\nabla\eta,\nabla\eta)\Gamma_{i}^{-1/2}(\Gamma_{i}-1)\ln^{2}(\Gamma_{i})\right)^{1/2}
\end{eqnarray*}
where in the last step the Cauchy-Schwarz inequality has been used for the nonnegative symmetric bilinear form $D^{2}f(\nabla u)$. 
Thus, after using Hölder's inequality, we find (recall (3.1))
\begin{eqnarray}
|(II)| &\le &  \left(\int_{\operatorname{spt}(\nabla\eta)\cap[\Gamma_{i}\ge2]}D^{2}f(\nabla u)(\nabla\partial_{i}u,\nabla\partial_{i}u)\Gamma_{i}^{-1/2}\eta^{2}\,dx\right)^{1/2} \nonumber\\[2ex]
&& \times\left(\int_{\operatorname{spt}(\nabla\eta)\cap[\Gamma_{i}\ge2]}D^{2}f(\nabla u)(\nabla\eta,\nabla\eta)\Gamma_{i}^{-1/2}(\Gamma_{i}-1)
\ln^{2}(\Gamma_{i})\,dx\right)^{1/2} \nonumber\\[2ex]
&\le & (\text{r.h.s. of } (3.2))
\end{eqnarray}
In a similar way we proceed with (I): we have
\begin{eqnarray*}
|\text{integrand of } (I)| &=& 2|\eta\partial_{i}u\Phi(\Gamma_{i})||D^{2}f(\nabla u)(\nabla\partial_{i}u,\nabla\eta)| \\[2ex]
&\le & c|\eta(\Gamma_{i}-1)^{1/2}\ln(\Gamma_{i})||D^{2}f(\nabla u)(\nabla\partial_{i}u,\nabla\eta)| \\[2ex]
&\le& c\left(D^{2}f(\nabla u)(\nabla\partial_{i}u,\nabla\partial_{i}u)\eta^{2}(\Gamma_{i}-1)\right)^{1/2}\\[2ex]
&& \times \left(D^{2}f(\nabla u)(\nabla\eta,\nabla\eta)\ln^{2}(\Gamma_{i})\right)^{1/2}
\end{eqnarray*}
and using Hölder's inequality again (recall (3.1))
\begin{eqnarray}
|(I)| &\le & c\left(\int_{\operatorname{spt}(\nabla\eta)\cap[\Gamma_{i}<2]}D^{2}f(\nabla u)(\nabla\partial_{i}u,\nabla\partial_{i}u)\eta^{2}(\Gamma_{i}-1)\,dx\right)^{1/2} \nonumber\\[2ex]
&& \times \left(\int_{\operatorname{spt}(\nabla\eta)\cap[\Gamma_{i}<2]}D^{2}f(\nabla u)(\nabla\eta,\nabla\eta)\ln^{2}(\Gamma_{i})\,dx\right)^{1/2} \nonumber\\[2ex]
&\le & (\text{r.h.s. of } (3.2))
\end{eqnarray}
We return to (3.5) noting that, on account of (3.6) it holds (recall (3.1))
\begin{equation}
(\text{l.h.s. of } (3.5))\ge c\int_{\Omega}D^{2}f(\nabla u)(\nabla\partial_{i}u,\nabla\partial_{i}u)\eta^{2}\lambda(\Gamma_{i})\,dx.
\end{equation}
Thus, combining (3.5) with (3.7)--(3.9), inequality (3.2) follows, which completes the proof of Lemma \ref{lem3.1}.
\end{proof}

\section{Proof of Theorem 1.1}

We want to prove Theorem \ref{thm1.1} in the setting of Remark 1.4, which means that we assume the conditions (1.2)--(1.5) for the density $f$ 
but replace (1.6) by ($p,q\in \mathbb{R}^n$, $|p|\le1$)
\begin{equation}
c|p|^{t-2}|q|^{2}\le D^{2}f(p)(q,q)\le C|p|^{\rho-2}|q|^{2}
\end{equation}
for exponents $2\le\rho\le t$, $t\le\rho+2$ and constants $c, C>0$.\\ 

Our arguments proceed as follows: let $i\in\{1,\dots,n\}$ and consider $\eta\in C_{0}^{1}(\mathbb{R}^{n})$ 
such that $\operatorname{spt}(\eta)\subset\overline{B}_{2R}$, $0\le\eta\le1$, $\eta=1$ 
on $\overline{B}_{R}$ and $|\nabla\eta|\le c/R$, where $B_{r}:=\{x\in\mathbb{R}^{n}:|x|<r\}$ for radii $r>0$. 
We will show that we can find a constant $c^{*}\in(0,\infty)$ for which
\begin{equation}
R^{-2}\int_{B_{2R}-B_{R}}|D^{2}f(\nabla u)|\Lambda(\Gamma)\,dx\le c^{*}, \quad 0<R<\infty,
\end{equation}
holds with $\Lambda$ defined in (3.1) and $\Gamma:=1+|\nabla u|^{2}$.\\ 

Accepting (4.2) for the moment, we quote (3.3) (recalling $\Lambda(\Gamma_{i})\le\Lambda(\Gamma)$, $\Gamma_{i}:=1+|\partial_{i}u|^{2}$) and get
\begin{equation*}
\int_{\mathbb{R}^{n}}D^{2}f(\nabla u)(\nabla\partial_{i}u,\nabla\partial_{i}u)\lambda(\Gamma_{i})\,dx<\infty.
\end{equation*}

Finally, we use inequality (3.2) with
\begin{equation}
\int_{B_{R}}D^{2}f(\nabla u)(\nabla\partial_{i}u,\nabla\partial_{i}u)\lambda(\Gamma_{i})\,dx =: \Theta_{i}(R)
\end{equation}
as a lower bound for the left-hand side of (3.2), whereas the right-hand side is bounded from above by the quantity ($c=\text{a constant independent of } R$)
\begin{equation*}
c\left(\int_{B_{2R}-B_{R}}D^{2}f(\nabla u)(\nabla\partial_{i}u,\nabla\partial_{i}u)\lambda(\Gamma_{i})\,dx\right)^{1/2}\cdot\left(R^{-2}\int_{B_{2R}-B_{R}}|D^{2}f(\nabla u)|\Lambda(\Gamma)\,dx\right)^{1/2}
\end{equation*}
which goes to zero as $R\rightarrow\infty$ on account of (4.2) and (4.3).\\ 

Thus we have $\lim_{R\rightarrow\infty}\Theta_{i}(R)=0$ and in conclusion
\begin{equation}
D^{2}f(\nabla u)(\nabla\partial_{i}u,\nabla\partial_{i}u)\lambda(\Gamma_{i})=0 \quad \text{on } \mathbb{R}^n
\end{equation}
for $i\in\{1,\dots,n\}$. The definition of $\lambda$ (see (3.1)) shows that $\lambda(\Gamma_{i}(x))=0$ if and only if $\partial_{i}u(x)=0$, 
and by Lemma 7.7, p. 152, in [6] it holds $\nabla\partial_{i}u=0$ almost everywhere on $[\partial_{i}u=0]$. 
However, on the set $[\partial_{i}u\ne0]$, (4.4) implies $\nabla\partial_{i}u=0$, hence $\nabla\partial_{i}u=0$ almost everywhere 
and by continuity everywhere on $\mathbb{R}^{n}$, being true for $i=1,\dots,n,$ completing the proof of Theorem \ref{thm1.1}.\\

Let us discuss our claim (4.2). On the set $[|\nabla u|\ge1]$ it holds by the definition of $\Lambda$ given in (3.1) and on account of (1.5)
\begin{align*}
|D^{2}f(\nabla u)|\Lambda(\Gamma) &= |D^{2}f(\nabla u)|\Gamma^{-1/2}(\Gamma-1)\ln^{2}(\Gamma) \\[2ex]
&\le |D^{2}f(\nabla u)||\nabla u|\ln^{2}(\Gamma) \\[2ex]
&\le c|\nabla u|^{1-\mu}|\nabla f(\nabla u)|.
\end{align*}
Note that in (1.5) we can evidently replace $\ln(1+|p|)$ by $\ln(1+|p|^{2})$ for $p\in\mathbb{R}^{n}$ such that $|p|\ge1.$ Using inequality (1.4), we get
\begin{equation}
|D^{2}f(\nabla u)|\Lambda(\Gamma)\le c\nabla u\cdot\nabla f(\nabla u) \quad \text{on } [|\nabla u|\ge1].
\end{equation}

On the set $[|\nabla u|<1]$ we observe using (3.1) and (4.1)
\begin{equation*}
|D^{2}f(\nabla u)|\Lambda(\Gamma)\le c|D^{2}f(\nabla u)|\ln^{2}(\Gamma)\le c|D^{2}f(\nabla u)||\nabla u|^{4}\le c|\nabla u|^{\rho+2}
\end{equation*}
since $\ln(1+\alpha)\le\alpha$ for $\alpha\ge0$. At the same time it holds (again for $[|\nabla u|<1]$)
\begin{equation*}
\nabla u\cdot\nabla f(\nabla u)=\int_{0}^{1}D^{2}f(r\nabla u)(\nabla u,\nabla u)\,dr \ge c\int_{0}^{1}r^{t-2}\,dr|\nabla u|^{t}=c|\nabla u|^{t},
\end{equation*}
and from our assumption $t\le\rho+2$ we find
\begin{equation}
|D^{2}f(\nabla u)|\Lambda(\Gamma)\le c|\nabla u|^{t} \quad \text{on } [|\nabla u|<1].
\end{equation}
Let $\mathcal{M}:\mathbb{R}^{n}\rightarrow[0,\infty)$,
\begin{equation}
\mathcal{M}(p):=\begin{cases}|\nabla f(p)||p|^{1-\mu},&|p|\ge1\\ |p|^{t},&|p|<1\end{cases}
\end{equation}
Then our calculations leading to (4.5) and (4.6) yield
\begin{equation}
|D^{2}f(\nabla u)|\Lambda(\Gamma)\le c\mathcal{M}(\nabla u)\le c\nabla u\cdot\nabla f(\nabla u),
\end{equation}
and (4.2) will follow as soon as we can show
\begin{equation}
R^{-2}\int_{B_{2R}-B_{R}}\mathcal{M}(\nabla u)\,dx\le c^{*}
\end{equation}
for a finite constant $c^{*}$ independent of $R$.\\ 

For proving (4.9) we select $\eta\in C_{0}^{1}(B_{\frac{5}{2}R}-B_{\frac{R}{2}})$, $0\le\eta\le1$ such that $\eta=1$ on $B_{2R}-B_{R},$ $|\nabla\eta|\le c/R$ 
and let $l\in\mathbb{N}$. Abbreviating $T_{R}:=B_{\frac{5}{2}R}-B_{\frac{R}{2}}$ we quote (4.8) and use equation (1.1) to get
\begin{eqnarray}\lefteqn{
\int_{T_{R}}\eta^{2l}\mathcal{M}(\nabla u)\,dx \le c\int_{T_{R}}\nabla u\cdot\nabla f(\nabla u)\eta^{2l}\,dx} \nonumber\\[2ex]
&=& -2lc\int_{T_{R}}\nabla f(\nabla u)\cdot\nabla\eta~\eta^{2l-1}u\,dx \nonumber\\[2ex]
&= & -2lc\left\{\int_{T_{R}\cap[|\nabla u|\ge1]}\nabla f(\nabla u)\cdot\nabla\eta~\eta^{2l-1}u\,dx + \int_{T_{R}\cap[|\nabla u|<1]}\dots \,dx\right\} \nonumber\\[2ex]
&= & -2lc((I)+(II)).
\end{eqnarray}

For $\epsilon, \delta>0$ we have (see (1.3))
\begin{align*}
|\text{integrand of } (I)| &\le |\nabla f(\nabla u)|\eta^{2l-1}|\nabla\eta||u| \\[2ex]
&= |\nabla f(\nabla u)|^{1-\epsilon}|\nabla f(\nabla u)|^{\epsilon}\eta^{2l-1}|\nabla\eta||u| \\[2ex]
&\le c|\nabla f(\nabla u)|^{1-\epsilon}|\nabla u|^{\epsilon(m-1)}\eta^{2l-1}|\nabla\eta||u| \\[2ex]
&\le \delta|\nabla f(\nabla u)||\nabla u|^{\frac{(m-1)\epsilon}{1-\epsilon}}\eta^{\frac{2l-1}{1-\epsilon}}
+c(\delta)|\nabla\eta|^{1/\epsilon}|u|^{1/\epsilon},
\end{align*}
and if we fix $\epsilon>0$ such that $\frac{\epsilon}{1-\epsilon}(m-1)\le1-\mu$ (recall $0\le\mu<1$), we obtain after increasing $l$ (recall (4.7)):
\begin{align*}
|\text{integrand of } (I)| &\le \delta|\nabla f(\nabla u)||\nabla u|^{1-\mu}\eta^{2l}+c(\delta)|\nabla\eta|^{1/\epsilon}|u|^{1/\epsilon} \\[2ex]
&\le \delta\mathcal{M}(\nabla u)\eta^{2l}+c(\delta)|\nabla\eta|^{1/\epsilon}|u|^{1/\epsilon}.
\end{align*}
Returning to (4.10) and choosing $\delta$ sufficiently small, it is shown
\begin{equation*}
\int_{T_{R}}\eta^{2l}\mathcal{M}(\nabla u)\,dx\le c\left\{|(II)|+\int_{T_{R}}|\nabla\eta|^{1/\epsilon}|u|^{1/\epsilon}\,dx\right\}.
\end{equation*}
The remaining term (II) is handled as follows (see (4.1) and (4.7)):
\begin{align}
|\text{integrand of } (II)| &\le |\nabla f(\nabla u)|\eta^{2l-1}|u||\nabla\eta| \nonumber\\[2ex]
&\le c|\nabla u|^{\rho-1}\eta^{2l-1}|\nabla\eta||u| \nonumber\\[2ex]
&\le \delta|\nabla u|^{t}\eta^{(2l-1)\frac{t}{\rho-1}}+c(\delta)|u|^{T}|\nabla\eta|^{T} \nonumber\\[2ex]
&\le \delta\mathcal{M}(\nabla u)\eta^{(2l-1)\frac{t}{\rho-1}}+c(\delta)|u|^{T}|\nabla\eta|^{T}
\end{align}
with exponent $T:=t/(t-\rho+1)$ as in Remark 1.4.\\ 

Finally, for $l$ sufficiently large, we can arrange that $(2l-1)t/(\rho-1)\ge2l$, and after appropriate choice of $\delta$ we deduce from (4.11)
\begin{equation}
\int_{T_{R}}\eta^{2l}\mathcal{M}(\nabla u)\,dx\le c\int_{T_{R}}\left(|u|^{T}|\nabla\eta|^{T}+|u|^{1/\epsilon}|\nabla\eta|^{1/\epsilon}\right)\,dx.
\end{equation}

We recall that $\epsilon>0$ just has to satisfy $\frac{\epsilon}{1-\epsilon}(m-1)\le1-\mu.$ 
Now we require in addition that $\epsilon\le1/T$ holds and use (1.7) on the right-hand side of (4.12) 
with the choice of $\beta=\beta(t,\rho)$ as defined in Remark 1.4. This gives (recall the choice of $\eta$) for $R\gg1$:
\begin{equation*}
\text{r.h.s. of } (4.12)\le cR^{n}\left[R^{\beta T-T}+R^{\frac{\beta}{\epsilon}-\frac{1}{\epsilon}}\right]\le cR^{n+\beta T-T}=cR^{2},
\end{equation*}
and (4.9) is a consequence of (4.12). This completes the proof of Theorem \ref{thm1.1} and of Remark 1.4, respectively.

\subsection*{Revisiting Example 2.1}

We finish by remarking that the arguments presented during the proof of Theorem \ref{thm1.1} are rather general, 
which means that they can easily be adjusted to particular cases, for instance it is possible to weaken condition (2.2) required in Example 2.1. 
We can state the following extension:

\begin{theorem} \label{thm4.1}
Let the density $f\in C^{2}(\mathbb{R}^{n})$ satisfy (1.2), (1.6) and (2.1), where in (2.1) we choose exponents $1<\tau\le r$ such that
\begin{equation}
r<\tau+1
\end{equation}
holds.\\ 

Then, if $u\in C^{2}(\mathbb{R}^{n})$ is an entire solution of (1.1) with the property (1.7), $\beta$ defined according to (1.8), it follows that $u$ is an affine function. Moreover, Corollary 1.2 and Corollary 1.3 extend to this situation.
\end{theorem}

\begin{proof}
We present the necessary adjustments. As before, it suffices to prove inequality (4.2). On the set $[|\nabla u|\ge1]$ we have by (2.1) (recall also (3.1))
\begin{align}
|D^{2}f(\nabla u)|\Lambda(\Gamma) &\le |D^{2}f(\nabla u)||\nabla u|\ln^{2}(\Gamma) \nonumber\\[2ex]
&\le c|\nabla u|^{r-2}|\nabla u|\ln^{2}(\Gamma)=c|\nabla u|^{r-1}\ln^{2}(\Gamma) \nonumber\\[2ex]
&\le c|\nabla u|^{\tau}
\end{align}
since
\begin{equation*}
\nabla u\cdot\nabla f(\nabla u)=\int_{0}^{1}D^{2}f(t\nabla u)(\nabla u,\nabla u)\,dt\ge c|\nabla u|^{\tau},
\end{equation*}
whereas (4.6) (now with $t=\rho$, see (1.6)) continues to hold. We replace (4.7) by
\begin{equation}
\tilde{\mathcal{M}}(p):=\begin{cases}|p|^{\tau},&|p|\ge1\\[1ex] |p|^{\rho},&|p|<1\end{cases}
\end{equation}
and deduce from (4.14) and (4.15)
\begin{equation}
|D^{2}f(\nabla u)|\Lambda(\Gamma)\le c\tilde{\mathcal{M}}(\nabla u)\le c\nabla u\cdot\nabla f(\nabla u).
\end{equation}
From (4.16) in turn it follows that our claim will be a consequence of (compare (4.9))
\begin{equation}
R^{-2}\int_{B_{2R}-B_{R}}\tilde{\mathcal{M}}(\nabla u)\,dx\le c^{*}
\end{equation}
again with a finite constant $c^*$ independent of $R$. For proving (4.17) we pass to inequality (4.10) discussing 
the quantity (II) exactly as before (recalling that now, for notational simplicity, we just assume that $\rho=t$). In addition 
we have for the integrand in (I) (compare (2.1) and (4.13))
\begin{align*}
|\nabla f(\nabla u)\cdot\nabla\eta~\eta^{2l-1}u| &\le |\nabla f(\nabla u)|\eta^{2l-1}|\nabla\eta||u| \\[2ex]
&\le c|\nabla u|^{r-1}\eta^{2l-1}|u||\nabla\eta| \\[2ex]
&\le \epsilon|\nabla u|^{\tau}\eta^{\frac{\tau}{r-1}(2l-1)}+c(\epsilon)|u|^{\frac{\tau}{\tau-r+1}}|\nabla\eta|^{\frac{\tau}{\tau-r+1}}
\end{align*}
After increasing $l$ and decreasing $\epsilon$ we get in place of (4.12)
\begin{equation}
\int_{T_{R}}\eta^{2l}\tilde{\mathcal{M}}(\nabla u)\,dx\le c\int_{T_{R}}\left(|u|^{a}|\nabla\eta|^{a}+|u|^{b}|\nabla\eta|^{b}\right)\,dx, \quad a:=\rho, \quad b:=\frac{\tau}{\tau-r+1} \, .
\end{equation}
Assume now that we have (1.7) with $\beta (\le1)$ from (1.8). Observe that in (2.1) we can increase $r$ such that $r$ is very close to $\tau+1$. 
Assuming this from the beginning, it follows that $a\beta-a\ge b\beta-b$ and (4.18) gives
\begin{equation*}
\int_{T_{R}}\eta^{2l}\tilde{\mathcal{M}}(\nabla u)\,dx\le cR^{n+\rho(\beta-1)}=cR^{2},
\end{equation*}
hence we have (4.17) and the proof is complete.
\end{proof}

\vspace{10mm}
\noindent
\begin{minipage}[t]{0.45\textwidth}
\textbf{Michael Bildhauer} \\
Saarland University \\
Department of Mathematics \\
P.O. Box 15 11 50 \\
66041 Saarbrücken, Germany \\
\texttt{bibi@math.uni-sb.de}
\end{minipage}
\hfill
\begin{minipage}[t]{0.45\textwidth}
\textbf{Martin Fuchs} \\
Saarland University \\
Department of Mathematics \\
P.O. Box 15 11 50 \\
66041 Saarbrücken, Germany \\
\texttt{fuchs@math.uni-sb.de}
\end{minipage}
\end{document}